\documentclass[10pt, a4paper]{article}

\usepackage{amsmath} 
\usepackage{amssymb} 
\usepackage{amsthm}
\usepackage{soul}
\usepackage{color}

\newtheorem{proposition}{Proposition}[section]
\newtheorem{rem}[proposition]{Remark}
\newtheorem{conjec}[proposition]{Conjecture}
\newtheorem{fait}[proposition]{Fact}
\newtheorem{definition}[proposition]{Definition}
\newtheorem{theo}[proposition]{Theorem}
\newtheorem{lem}[proposition]{Lemma}
\newtheorem{cor}[proposition]{Corollary}

\newcommand{\ad}{\operatorname{ad}}
\newcommand{\id}{\operatorname{Id}}
\newcommand{\rk}{\operatorname{rk}}

\newcommand{\gk}{\mathfrak{g}}
\newcommand{\hk}{\mathfrak{h}}
\newcommand{\bk}{\mathfrak{b}}
\newcommand{\ak}{\mathfrak{a}}
\newcommand{\ck}{\mathfrak{c}}
\newcommand{\uk}{\mathfrak{u}}
\newcommand{\nk}{\mathfrak{n}}
\newcommand{\mk}{\mathfrak{m}}
\newcommand{\dk}{\mathfrak{d}}
\newcommand{\ik}{\mathfrak{i}}
\newcommand{\jk}{\mathfrak{j}}
\newcommand{\K}{\mathbb{K}}
\newcommand{\N}{\mathbb{N}}
\newcommand{\ov}{\overline}

\title{Cartan subrings in ranked soluble Lie rings}
\author{Jules Tindzogho Ntsiri and Samuel Zamour}

\begin{document}

\maketitle

\abstract{We prove the existence of Cartan subrings in soluble ranked Lie rings.}

\section{Introduction}

More than forty years after Rosengarten's PhD thesis \cite{Ros}, ranked Lie rings seem to have, in recent years, attracted renewed interest among model theorists. This renewed interest is due to the work of Deloro and Tindzogho in \cite{DT1}. Notably, because they took-up and further Rosengarten's results by adopting a more algebraic conceptualization closer to that of non-logicians, but also because they propose a list of questions for research, some of which have captured the attention of researchers. Moreover, in \cite{DT2}, the two authors provided an answer to question 8 on this list and gave an analog to the Lie-Kolchin theorem for ranked Lie rings, which will be essential for our study.

This paper is focused on question 10, and what is presented here can be summarized in two results. The first one provides the existence of \textit{Cartan subrings}, i.e., definable connected nilpotent self-normalizing subrings, in our context:\\

\noindent{\bf Corollary \ref{Cartan subrings exist}.} 
{\it A connected non-nilpotent soluble ranked Lie ring $\gk$ such that $\gk'$ is nilpotent has Cartan subrings.}\\

The second result gives a characterization of Cartan subrings:\\

\noindent{\bf Theorem \ref{main theorem}.}
{\it Let $\gk$ be a connected soluble ranked Lie ring such that $\operatorname{char}(\gk)>\rk(\gk)$ and let $\ck\sqsubseteq\gk$ be a subring. The following items are equivalent:
\begin{itemize}
    \item [(i)] $\ck$ is a minimal def-abnormal subring;
    \item [(ii)] $\ck$ is a Cartan subring;
    \item [(iii)] $\ck$ is an Engel-minimal subring.
\end{itemize}}

Our approach will consist of drawing on Frécon's work on \emph{abnormal subgroups} in soluble groups of finite Morley rank \cite{Fre} (we also refer to \cite[sec.8, chap. I]{ABC}, which gives a concise and elegant treatment of Frécon's results). 

\begin{definition}
 Let $G$ be a group and $H$ a subgroup of $G$. We say that $H$ is abnormal in $G$ if, for all $g\in G$, $g$ belongs to $\langle H, H^g\rangle$. 
\end{definition}
 Using abnormal subgroups, Frécon proved the existence and conjugacy of \emph{Carter subgroups} in soluble groups of finite Morley rank. Carter subgroups are those that are nilpotent and self-normalizing. We would like to point out that the above definition does not apply to our context, because it relies on the use of inner automorphisms which are meaningless for Lie rings. Instead, we adopt the definition of abnormality provided by \cite{stitz} (see Definition \ref{definition of abnormality}). Notice that Barnes proved that, for finite dimensional Lie algebras, Cartan subalgebras are in fact abnormal. Moreover, abnormality behaves well with respect to quotients: the image of an abnormal subring is abnormal, and we can lift an abnormal subring from a quotient to an abnormal subring (see Lemma \ref{abnormal and quotient}). This is crucial when reasoning by induction, and a priori not guaranteed when starting from the definition of a Cartan subring. Nevertheless, adopting a more direct treatment in proving Corollary \ref{Cartan subrings exist}, based on the original definition of Cartan subrings, seems a desirable goal for a future study. 

In \cite{Bar}, Barnes also proved that Cartan subalgebras are minimal Engel subalgebras (see \cite[Theorem 1]{Bar}). Through Theorem \ref{main theorem}, we aim to establish a similar result. To achieve this, Fact \ref{Frattini argument} plays a very important role in obtaining the nilpotency of the Frattini ideal, which is essential for the treatment of Engel subrings.

Now let us say a few words about the assumptions. The main assumption in this article, is certainly the nilpotency of the derived subring of a soluble ranked Lie ring. This assumption, assumed early in the paper  (Proposition \ref{Cartan in 2-solvalble}), conditions a large part of this article.  Indeed, the nilpotency of the derived subring is not guaranteed as in the case for groups. In \cite{DT2}, it is obtained under an assumption on the characteristic (Fact \ref{g' is nilpotent}). This fact will be used in the last part of this work, from Lemma \ref{Fitting is the set of ad-nilpotent} to the end.

We want to emphasize that our work has no connection with the study of the Cherlin-Zilber Conjecture, which states that an infinite simple group of finite Morley rank is algebraic over an algebraically closed field. However, it should contribute to answering the "logged" version:

\begin{conjec}\cite{DT1}
Let $\gk$ be an infinite simple Lie ring of finite Morley rank. Suppose the characteristic is sufficiently large. Then $\gk$ is a simple Lie algebra over an algebraically closed field.    
\end{conjec} 

Furthermore, it is difficult to see how to establish, as is the case between Lie groups and Lie algebras, a Lie-Chevalley correspondence between ranked groups and ranked Lie rings (see \cite[Question 14]{DT1}).

In the second section, where we give the main definitions and notation used, we recall some linearization results and give other tools. In Section 3, under the assumption of nilpotency of the derived subring of a connected soluble ranked Lie ring, we develop the theory of abnormal subrings, prove their existence, and show that minimal ones are Cartan subrings. Section 4 is brief and focuses on the Frattini subring. Finally, in Section 5, we deal with Engel subrings in order to give another characterization of Cartan subrings.

\section{Definitions and notations}
\begin{itemize}
\item A Lie ring is an abelian group $\gk$ equipped with a bi-additive and anti-symmetric function, $[\cdot,\cdot]:\gk\times\gk\to\gk$ that satisfies the Jacobi identity:
\[[a,[b,c]]=[[a,b],c]+[b,[a,c]]\]
\item A Lie ring can be seen as a Lie algebra devoid of its vector space structure. Therefore, the notions of subring, ideal, abelianity, solubility and nilpotency keep the same definition.
\item A ranked Lie ring is the one which is definable in a finite Morley rank structure. Here, the Morley rank is a dimension defined on definable sets. We refer the reader to \cite{BN} (especially chapter 4) for more details on the matter. 
\item $\rk(\gk)$ will denote the Morley rank of $\gk$.
\item The characteristic of a Lie ring is its exponent as an abelian group.
\item As in \cite{DT1,DT2}, we use $\leqslant$ for subgroups, $\sqsubseteq$ for subrings, $\lhd$ for ideals and $(+)$ means with finite intersection.
\end{itemize}

Throughout the following, $\gk$ stands for a ranked Lie ring, and by \emph{subring} we mean \emph{Lie subring}.
A \emph{definable $\gk$-module} is a definable abelian group $V$ on which $\gk$ acts definably as a Lie ring. 

\begin{definition}
A $\gk$-module $V$ is said to be \emph{$\gk$-minimal} if it is definable connected and it does not have an infinite definable $\gk$-submodule. It is said to be \emph{irreducible} if, in addition, it has no finite $\gk$-submodule.
\end{definition}

\begin{fait}(see \cite[Lemma 4.10]{ABC})\label{irreducibility}
Let $\gk$ be an abelian ranked Lie ring and $V$ be a $\gk$-minimal module. Suppose that the action is non-trivial. Then $V$ is $\gk$-irreducible.     
\end{fait}

\begin{fait}\label{Lie theorem}\cite[Corollary A1]{DT2}
Let $\gk$ be a connected soluble ranked Lie ring. Let $V$ be a faithful $\gk$-minimal module. Suppose $\operatorname{char} V > \rk( V)$. Then $\gk$ is abelian and the configuration is definably linear: there is an infinite definable field $\K$ such
that $V \cong \K_+$ and $\gk\hookrightarrow \K\id_V $, all definably.
\end{fait}

\begin{lem}\label{irreducibility of minimal modules}
    Let $\gk$ be a connected soluble ranked Lie ring such that $\operatorname{char}(\gk)> \rk(\gk)$ and $V$ a $\gk$-minimal module on which $\gk$ acts non-trivially. Then $V$ is $\gk$-irreducible.
\end{lem}  

\begin{proof}
 This follows immediately from Fact \ref{Lie theorem} above.
\end{proof}

\begin{fait}\cite[Corollary A2]{DT2}\label{g' is nilpotent}
Let $\gk$ be a connected soluble ranked Lie ring. Suppose $(\gk;+)$ has prime exponent $> \rk(\gk)$.
\begin{itemize}
\item [(i)] There is a largest definable, connected, nilpotent ideal $F^{\circ}(\gk) \lhd\gk$.
\item [(ii)] $\gk'$ is nilpotent.
\item [(iii)] If $\uk \sqsubseteq \gk$ is a definable, connected subring consisting of ad-nilpotent elements (viz., $(\forall x\in \uk)(\exists n\in\N)(\ad_x^n(\gk)=0$), then $\uk\sqsubseteq F^{\circ}(\gk)$.
\end{itemize}
\end{fait}

\begin{definition}
Let $\gk$ be a ranked Lie ring. A Lie subring $\ck\sqsubseteq\gk$ is called a \emph{Cartan subring} of $\gk$ if it is definable, connected, nilpotent and self-normalizing.
\end{definition}

\section{(Def)-Abnormal subrings}

In the following, we say that two subrings $\hk_1, \hk_2$ of a Lie ring $\gk$ are \emph{conjugated}, if there exists $x\in \gk$ such that $\ad_x$ is nilpotent and induces an  \emph{inner} automorphism $\exp(\ad_x)$ of $\gk$ with $\hk_2=\exp(\ad_x)(\hk_1)$. For example, if $\ik\lhd \gk$ is an abelian ideal, then for any $x\in \ik$, $\exp(\ad_x)=\id+\ad_x$. Here, the exponential is defined in an usual manner without any assumption on the characteristic.

In the result below, we prove the existence and conjugacy of Cartan subrings in connected $2$-soluble ranked Lie rings with a minimal derived subring.

\begin{proposition}\label{Cartan in 2-solvalble}
Let $\gk$ be a connected non-nilpotent soluble ranked Lie ring such that $\gk'$ is $\gk$-minimal. Let $a\in \gk\setminus{C_{\gk}(\gk')}$ and $\ak=C_{\gk}^{\circ}(a)$. Then
\begin{itemize}
    \item[(i)] $\gk=\gk'(+)\ak$;
    \item[(ii)] $\ak$ is a Cartan subring of $\gk$;
    \item[(iii)] if $\uk$ is a subring containing $\ak$, then $\uk=\ak$ or $\uk=\gk$;
    \item[(iv)] if $\bk\sqsubseteq\gk$ is definable, connected and such that $\gk=\gk'(+) \bk$, then $\bk$ is $\exp(\gk')$-conjugated to $\ak$. In particular, $\bk$ is a Cartan subring of $\gk$.
\end{itemize}  
\end{proposition}

\begin{proof}
\begin{itemize}
\item [(i)] First note that $\gk''=0$ by $\gk$-minimality of $\gk'$, so $\gk'$ is abelian. This implies that $(\ak\cap \gk')^{\circ}$ is an ideal of $\gk$,  and by $\gk$-minimality we have $(\ak\cap\gk')^{\circ}=0$. In particular, the restriction of $\ad_a$ to $\gk'$ has a finite kernel, hence surjective. Then $[a,\gk]=[a,\gk']=\gk'$ and this proves (i). 
\item [(ii)] Since $(\ak\cap \gk')^{\circ}=0$, $\ak$ is abelian hence nilpotent. Moreover, since $\gk'$ is $\ak$-minimal and $\ak$ is abelian acting non-trivially ($\gk$ is not nilpotent), by Fact \ref{irreducibility}, $\gk'$ is $\ak$-irreducible. In particular, we have $\gk'\cap \ak=0$. Thus $C_{\gk}(\ak)=(C_{\gk}(\ak)\cap\gk')\oplus\ak$ and so $C_{\gk}(\ak)=\ak$ since $\gk$ is not nilpotent. It remains to show that $\ak=N_{\gk}(\ak)$. Let $x\in N_{\gk}(\ak)$, we have $[x,\ak]\leqslant(\ak\cap\gk')=0$, so $N_{\gk}(\ak)\leqslant C_{\gk}(\ak)$ and we are done. 
\item[(iii)] Let $\uk$ be a subring containing $\ak$; then $\uk=(\uk\cap \gk')(+)\ak$. Since $\gk'$ is $\ak$-irreducible and $\ak\subseteq N_{\gk}(\uk)$, we have two cases : $\uk=\gk$ or $\uk=\ak$.
\item [(iv)] $\bk$ does not centralize $\gk'$ otherwise $\gk$ is nilpotent. Let $b\in\bk\setminus{C_{\gk}(\gk')}$; so $b=y+a_1$ with $y\in \gk'$ and $a_1\in\ak$. As above, $\ad_{a_1}:\gk'\to\gk'$ is surjective, so there exists $z\in \gk'$ such that $y=[a_1,z]$ so that $\exp(\ad_{z})(a_1)=a_1+[z,a_1]=b$. But $\exp(\ad_{z})$ is a Lie automorphism of $\gk$, and thus $\ck=\exp(\ad_{z})(\ak)$ is an abelian Cartan subring containing $b$. It follows that $\ck\sqsubseteq C_{\gk}^{\circ}(b)=\bk$ and the equality holds by rank consideration.
\end{itemize}
\end{proof}

This theorem generalizes \cite[Theorem 2]{DT1} to any connected $2$-soluble ranked Lie ring with a minimal derived subring. To prove the existence of Cartan subring in any soluble ranked Lie ring we need to introduce abnormal subrings.

\begin{definition}\label{definition of abnormality}
 Let $\gk$ be a (ranked) Lie ring. A (definable) subring $\ak$ is called \emph{(def)-abnormal} if for all (definable) subrings $\mathfrak{u}$ such that $\ak\sqsubseteq\mathfrak{u}$, we have $N_{\gk}(\mathfrak{u})=\mathfrak{u}$.
\end{definition}

\begin{rem}\label{Cartan is abnormal}
 Cartan subrings involved in Proposition \ref{Cartan in 2-solvalble} are def-abnormal, and this mention is important in the following.
\end{rem}

\begin{lem}\label{abnormal and quotient}
Let $\gk$ be a (ranked) Lie ring and $\mathfrak{i}$ a (definable) ideal. If $\ck$ is (def-)abnormal, then $\ck+\ik/\ik $ is (def)-abnormal. Conversely,  if $\overline{\ck}$ is (def-)abnormal in $\gk/\ik$, then its  preimage $\mathfrak{c}=\pi^{-1}(\overline{\mathfrak{c}})$ is (def)-abnormal in $\gk$.
\end{lem}

\begin{proof}
Suppose that $\ck$ is (def-)abnormal. Let $\uk$ be a (definable) subring such that $\ck+\ik\sqsubseteq\uk+\ik$. If $x+\ik/\ik$ normalizes $\uk+\ik/\ik$, then $x$ normalizes $\uk+\ik$, hence $x$ lies in $\uk+\ik$ since $\ck$ is (def)-abnormal. Thus $x+\ik/\ik\in \uk +\ik/\ik$ and $N_{\gk/\ik}(\uk+\ik)=\uk+\ik/\ik$ which means that $\ck+\ik/\ik$ is (def-)abnormal. 

Conversely, if $x$ normalizes $\uk$ a (definable) subring containing $\ck$, then $x+\ik\in N_{\gk/\ik}(\uk/\ik)=\uk/\ik$ and $x\in \uk$.  We can conclude.
\end{proof}
Now we prove that def-abnormal subrings are connected.
\begin{proposition}\label{abnormal is def and connected}
Let $\gk$ be a soluble connected ranked Lie ring such that $\gk'$ is nilpotent and let $\ck$ be a def-abnormal subring. Then $\ck$ is connected.
\end{proposition}

\begin{proof}
We proceed by induction on $\rk(\gk)$. Let $\ik$ be a definable connected minimal ideal contained in $\gk'$. By induction $\ck+\ik/\ik$ is definable and connected, and so is $\ck+\ik$. We can assume that $\gk=\ck+\ik$. In particular, $\ck\cap \ik$ is a finite ideal by minimality of $\ik$. Since $\gk'$ is nilpotent, $\gk'\leqslant C_{\gk}(\ik)$. We can also assume that $\ik$ is not central otherwise $\ik\sqsubseteq \ck$ and $\gk=\ck$. Therefore $\gk/\gk'$ acts non-trivially on $\ik$; by Fact \ref{irreducibility}, $\ik$ is $\gk/\gk'$-irreducible. Thus we have $\ck\cap \ik=0$ and we are done since $\gk=\ik\oplus \ck$. 
\end{proof}

\begin{fait} \label{nilpotency criterion}(Hall's criterion for nilpotency \cite[Theorem 2]{Ch})
Let $\gk$ be a Lie ring and $\ik$ an ideal. Then $\gk$ is nilpotent if and only if $\ik$ and $\gk/\ik'$ are nilpotent.
\end{fait}

\begin{proposition}\label{g-irreducible quotient}(compare with \cite[Lemma 8.15, chap. I]{ABC})
Let $\gk$ be a connected, soluble and non-nilpotent ranked Lie ring such that $\gk'$ is nilpotent. Then $\gk$ has a definable centerless quotient $\overline{\gk}$ such that $\overline{\gk}'$ is abelian and $\gk$-irreducible.   
\end{proposition}

\begin{proof}
    Let $\ik$ be a connected definable ideal such that $\gk/\ik$ is non-nilpotent and maximal as such. Replacing $\gk$ with $\gk/\ik$, we can assume that any non-trivial quotient of $\gk$ by a connected definable ideal is nilpotent; thus $Z(\gk)$ is finite and $\gk/Z(\gk)$ is centerless. Replacing $\gk$ with $\gk/Z(\gk)$ we can also assume that $\gk$ is centerless. 
    
    If $\gk'$ is not abelian, then  $\gk/\gk''$ is nilpotent, and so is $\gk$ by Hall's nilpotency criterion (Fact \ref{nilpotency criterion}); a contradiction. So $\gk'$ is abelian.

    Considering the DDC on the descending central series, let $k$ be minimal such that $\gk^k=\gk^{k+1}$. Since $\gk$ is non-nilpotent, $\gk^{k}\neq0$. We claim that $\gk^k$ is $\gk$-irreducible. Let $\ik\subseteq \gk^k$ be an ideal. If $\ik$ is infinite, $\gk/\ik^{\circ}$ is nilpotent and $\gk^k\subseteq \ik$; the equality holds. If $\ik$ is finite, it is central hence trivial. Therefore $\gk^k$ is $\gk$-irreducible.

    Now we prove that $\gk'=\gk^k$ to finish the proof. Towards a contradiction, suppose the opposite. Consider $\overline{\gk}=\gk/\gk^{k}$ and let $\overline{\ak}\leqslant \overline{\gk'}$ be $\gk$-minimal with $\ak>\gk^k$ as its preimage. Since $\overline{\gk}$ is nilpotent, $\overline{\ak}$ is central. In particular, $[\ak,\gk]\leqslant \gk^k$. Thus, for any $x\in\gk$, the centralizer $C_{\ak}(x)$ is an ideal since $\gk'$ is abelian and $\ak\subseteq \gk'$. It follows that $[\gk,C_{\ak}(x)]$ is a connected definable ideal contained in $\gk^k\cap C_{\ak}(x)$, and properly contained in $\gk^{k}$ if one chooses $x\in\gk-C_{\gk}(\gk^k)$. By $\gk$-irreducibility of $\gk^k$ we get $[\gk,C_{\ak}(x)]=0$, that is $C_{\ak}(x)\leqslant Z(\gk)=0$. This implies that $\ad_x:\ak\to\gk^k$ is  injective; a contradiction as $\gk^k<\ak$.
\end{proof}

\begin{proposition} \label{abnormal exits}
    Let $\gk$ be a connected non-nilpotent and soluble  ranked Lie ring such that $\gk'$ is nilpotent. Then $\gk$ has a proper def-abnormal subring. 
    \end{proposition}

    \begin{proof}
This follows from Proposition \ref{g-irreducible quotient} above combined with Lemma \ref{abnormal and quotient} and Remark \ref{Cartan is abnormal}.
    \end{proof}

\begin{proposition} \label{abnormality criterion}
Let $\gk$ be a connected soluble ranked Lie ring such that $\gk'$ is nilpotent, $\ik$ be a connected minimal ideal, and $\hk$ be a definable subring of $\gk$ containing  $Z^{\circ}(\gk)$, such that $\gk = C_{\gk}(\ik)+\hk$. Then $\hk$ is def-abnormal in $\hk+\ik$. 
\end{proposition}

\begin{proof}
We can assume that $\ik$ is not contained in $\hk$. Notice that $\ik$ is $\hk$-minimal. Also, $\gk'$ centralizes $\ik$ since it is nilpotent and $\ik$ is not containing in $Z^{\circ}(\gk)$. Thus $\gk/\gk'$ acts non-trivially on $\ik$ and by Fact \ref{irreducibility}, $\ik$ is $\gk/\gk'$-irreducible. It follows that $\ik\cap Z(\gk)=\ik\cap \hk=0$.

Now let $\uk\sqsubseteq \ik+\hk$ be a definable subring containing $\hk$. Set $\bk=\ik\oplus\hk$ and $\nk=N_{\bk}(\uk)$. We have $\nk=(\nk\cap \ik)\oplus\hk$ and $\nk\cap \ik$ is $\hk$-invariant and even $\gk$-invariant. If $\nk\cap \ik$ is finite, it is central in $\gk$, hence trivial and we get $\nk=\hk=\uk$ so we are done. Therefore, by $\hk$-minimality of $\ik$ we can assume that $\ik$ normalizes $\uk$. In particular we have $\nk=\bk$. Similarly, if $\ik\cap \uk$ is finite, it is trivial and we get $\hk=\uk$, so we are done. Thus, by $\hk$-minimality we have $\ik\sqsubseteq \uk$ so $\nk=\bk=\uk$.
\end{proof}

\begin{proposition}\label{abnormality is transitive}
Let $\gk$ be a connected soluble ranked Lie ring such that $\gk'$ is nilpotent, $\ak$ be a def-abnormal subring of $\gk$, and $\hk$ be a subring of $\ak$ which is 
def-abnormal in $\ak$. Then $\hk$ is def-abnormal in $\gk$. 
\end{proposition}

\begin{proof}
Let $\gk$ be a counterexample of minimal rank with $\hk$ of maximal rank. Since $\ak+\gk'/\gk'$ and $\hk+\gk'/\gk'$ are def-abnormal and $\gk/\gk'$ is abelian, we have $\gk=\gk'+\ak$ and $\ak+\mathfrak{g'}=\gk'+\hk$, hence $\gk=\gk'+\hk$. Let $\ik$ be a connected minimal ideal contained in $\gk'$. Then $\gk'\subseteq C_{\gk}(\ik)$ ($\gk'$ is nilpotent), and we get $\gk=C_{\gk}(\ik)+\hk$. But $\ak$ is def-abnormal in $\gk$, so $Z(\gk)\subseteq\ak$, hence $Z(\gk)\subseteq \hk$. Therefore, Proposition \ref{abnormality criterion} implies that $\hk$ is def-abnormal in $\hk+\ik$. 

In addition, $\hk+\ik/\ik$ is def-abnormal in $\gk/\ik$ (by induction), so $\hk+\ik$ is def-abnormal in $\gk$ by Lemme \ref{abnormal and quotient}. Now, towards a contradiction, we show that $\hk$ is def-abnormal in $\gk$.

Let $\uk$ be a definable subring of $\gk$ that contains $\hk$ and let $x\in N_{\gk}(\uk)$. In particular $x$ normalizes $\ik+\uk$, and then lies in $\ik+\uk$ since $\ik+\hk$ is def-abnormal in $\gk$. Write $x=i+u$ with $i\in \ik$ and $u\in \uk$. But $i$ normalizes $\uk$, so $i\in N_{\hk+\ik}(\uk\cap(\hk+\ik))$. By def-abnormality of $\hk$ in $\hk+\ik$, we have $i\in \uk\cap(\hk+\ik)$, and we get $x\in \uk$, thus $N_{\gk}(\uk)=\uk$ and the conclusion follows.
\end{proof}

\begin{proposition}\label{minimal abnorma is Cartan}
Let $\gk$ be a connected non-nilpotent soluble ranked Lie ring such that $\gk'$ is nilpotent. Then any minimal def-abnormal subring is a Cartan subring. 
\end{proposition}

\begin{proof}
Let $\ak$ be a minimal def-abnormal subring (such a subring exists by Proposition \ref{abnormal exits}). It is immediate that $\ak$ is self-normalizing. It remains to show that $\ak$ is nilpotent. Towards a contradiction, suppose not. Then by Proposition \ref{abnormal exits} $\ak$ has a proper def-abnormal subring. The latter is def-abnormal in $\gk$ by Proposition \ref{abnormality is transitive}, which contradicts the minimality of $\ak$.
\end{proof}

\begin{cor}\label{Cartan subrings exist}
A connected non-nilpotent soluble ranked Lie ring $\gk$ such that $\gk'$ is nilpotent has Cartan subrings.
\end{cor}

The following deals with the conjugacy; in a Lie algebra $\gk$, for a subalgebra $\uk$, we denote by $\mathcal{I}(\gk,\uk)$ the subgroup of automorphisms generated by $\exp(\ad_x)$ for all $x\in \uk$ for which the exponential makes senses.

\begin{fait}\cite[Lemma 5]{Bar}\label{inner automorphisms}
Let $\gk$ be a soluble Lie algebra over a field $F$. If $\operatorname{char} F=p \neq0$, suppose further
that $\gk'$ is $(p-1)$-Engel, i.e., $\ad_x^{p-1}=0$ on $\gk$ for all $x\in \gk'$. Then
\begin{itemize}
\item [(i)] for all $x\in\gk', \exp \ad_x$ exists,
\item [(ii)] If $\uk\sqsubseteq \gk$, every $\alpha\in \mathcal{I}(\uk, \uk')$ has an extension $\alpha^*\in\mathcal{I}(\gk, \uk')$,
\item [(iii)] if $\ik\lhd\gk$, every $\beta\in\mathcal{I}(\gk/\ik, (\gk/\ik)')$ is induced by some $\beta^*\in\mathcal{I}(\gk, \gk')$.
\end{itemize}
\end{fait}

\begin{proposition}\label{conjugacy of minimal abnormal}
Let $\gk$ be a connected non-nilpotent soluble ranked Lie ring of characteristic $p$ and such that $\gk'$ is $(p-1)$-Engel. Then the nilpotent def-abnormal subrings  of $\gk$ are conjugated under $\mathcal{I}(\gk,\gk')$.    
\end{proposition}

\begin{proof}
We proceed by induction on $\rk(\gk)$. Let $\ak$ and $\bk$ be any two def-abnormal subrings of $\gk$. If $Z(\gk)$ is infinite, then by applying induction to $\gk/Z(\gk)$, $\ak/Z(\gk)$ and $\bk/Z(\gk)$ are conjugated under $\mathcal{I}(\gk/Z(\gk),(\gk/Z(\gk))')$, so by Fact \ref{inner automorphisms} (iii) $\ak$ and $\bk$ are conjugated under $\mathcal{I}(\gk,\gk')$. Now suppose that $Z(\gk)$ is finite, so replacing $\gk$ with $\gk/Z(\gk)$ we can assume that $Z(\gk)=0$.

Let $\ik$ be a $\gk$-minimal ideal. By induction, $\ak+\ik/\ik$ and $\bk+\ik/\ik$ are conjugated under $\mathcal{I}(\gk/\ik,(\gk/\ik)')$. Again by Fact \ref{inner automorphisms} (iii) $\ak+\ik$ and $\bk+\ik$ are conjugated under $\mathcal{I}(\gk,\gk')$. Since the image of a def-abnormal subring by an (inner) automorphism of $\gk$ is still def-abnormal, we can suppose that $\gk=\ak+\ik=\bk+\ik$. But $C_{\ak}(\ik)$ is an ideal of $\ak$ and the latter is nilpotent, so $C_{\ak}(\ik)=0$ otherwise we have $0\neq Z(\ak)\cap C_{\ak}(\ik)\subseteq Z(\gk)$. It follows that $C_{\gk}(\ik)=\ik$ and since $\gk'$ centralises $\ik$, we have $\gk'=\ik$ and we can conclude by Proposition \ref{Cartan in 2-solvalble}.
\end{proof}

\section{Frattini subring}

\begin{definition}
The \emph{Frattini} subring $\Phi(\gk)$ of a connected ranked Lie ring $\gk$ is the intersection of all the proper definable connected maximal subrings.
\end{definition}
Notice that contrary to the analogous notion in group theory, it is not obvious that the Frattini subring is in fact an ideal.
\begin{fait}\label{Frattini fact}
Let $\gk$ be a connected ranked Lie ring, and $\hk\sqsubseteq \gk$ a definable subring such that $\gk=\Phi(\gk)+\hk$. Then $\gk=\hk$.
\end{fait}
The following proposition shows that the Frattini subring of a connected soluble ranked Lie ring is an ideal; we adapt \cite[Lemma 3.4]{BG} to our context.
\begin{proposition}\label{Frattini is an ideal}
Let $\gk$ be a connected soluble ranked Lie ring such that $\gk'$ is nilpotent. Then $\Phi(\gk)$ is an ideal of $\gk$.
\end{proposition}

\begin{proof}
We  proceed by induction on $\rk(\gk)$. The result is trivial if $\gk$ is abelian (or $\rk(\gk)=1$). Let $\ik$ be a $\gk$-minimal ideal of $\gk$. Denote by $\Phi_{\ik}$ the intersection of all the proper definable connected maximal subrings containing $\ik$. By induction, we find $\Phi_{\ik}/\ik=\Phi(\gk/\ik)$ is an ideal of $\gk/\ik$ so $\Phi_{\ik}\lhd \gk$.

If $\ik\subseteq\Phi(\gk)$, then $\Phi(\gk)=\Phi_{\ik}\lhd\gk$. Now suppose that $\ik\not\subseteq\Phi(\gk)$ and let $\mk$ be a definable connected maximal proper subring not containing $\ik$. Then $\gk=\ik+\mk=C_{\gk}(\ik)+\mk$. In particular $\ik$ is $\mk$-minimal and  $C_{\gk}(\ik)=\ik+C_{\mk}(\ik)$. 

If $\ik<C_{\gk}^{\circ}(\ik)$ then $C_{\mk}(\ik)$ is an infinite ideal of $\gk$. It follows that every definable connected proper maximal subring contains an infinite ideal $\jk$. This implies that $\Phi(\gk)=\bigcap\limits_\jk \Phi_{\jk}\lhd \gk$.

Suppose $\ik=C_{\gk}^{\circ}(\ik)$. Then $C_{\mk}(\ik)$ is a finite ideal. Note that $\ik\cap\mk=0$ : by minimality of $\ik$ and nilpotency of $\gk'$, $\gk'\subseteq C_{\gk}^{\circ}(\ik)=\ik$; thus, $\mk'\subseteq \ik\cap \mk$ and $\mk'$ is trivial by connectedness, so we can use Fact \ref{irreducibility}.  
We claim that $\Phi(\gk)$ is contained in $C_{\mk}(\ik)$. Indeed, let $x\in\mk\setminus C_{\mk}(\ik)$; there is  $a\in \ik$ such that $[a,x]\neq0$. In particular, $\ad_a$ is $2$-nilpotent. Thus $\exp(\ad_a)=\id +\ad_a$ is an automorphism of $\gk$ and $\mk_1=\exp(\ad_a)(\mk)$ is a definable connected proper maximal subring. If $x\in\mk_1\cap\mk$, then $x=y+[a,y]$ with $y\in \mk$, so $[a,y]\in \ik\cap\mk$. But $\ik\cap \mk=0$, so $[a,y]=0$. This implies $x=y$ and $[a,x]=0$, a contradiction. Therefore $x\notin \Phi(\gk)$ and the claim is proved. Since $C_{\mk}(\ik)$ is a finite ideal, it is central, so $\Phi(\gk)\lhd \gk$.
\end{proof}

\begin{fait}\cite{Zam2}\label{Frattini argument}{\em (Frattini argument)}

Let $\gk$ be a connected ranked Lie ring, $\ck$ a connected Cartan subring of a definable connected ideal $\ik$. Then $\gk=\ik+N_{\gk}(\ck)$.
\end{fait}

\begin{proposition}\label{frattini is nilpotent}
Let $\gk$ be a connected soluble ranked Lie ring such that $\gk'$ is nilpotent. Then $\Phi^{\circ}(\gk)$ is nilpotent.
\end{proposition}

\begin{proof}
Towards a contradiction, suppose that $\Phi^{\circ}(\gk)$ is not nilpotent. By Proposition \ref{Cartan subrings exist} $\Phi^{\circ}(\gk)$ has connected Cartan subrings. Since $\Phi^{\circ}(\gk)$ is a definable connected ideal (Proposition \ref{Frattini is an ideal}), Fact \ref{Frattini argument} implies that $\gk=\Phi^{\circ}(\gk)+N_{\gk}(\ck)$, with $\ck$ a connected Cartan subring of $\Phi^{\circ}(\gk)$. By Fact \ref{Frattini fact}, it follows that $\gk=N_{\gk}(\ck)$ and we have $\ck=\Phi^{\circ}(\gk)$, that is $\Phi^{\circ}(\gk)$ is nilpotent, a contradiction.
\end{proof}
We define the \emph{Fitting} ideal of $\gk$ to be the largest  nilpotent ideal (it is denoted by $F(\gk)$). It exists and it is definable by \cite{Zam1}. We obtain a characterization of the Fitting ideal modulo the Frattini ideal (compare with \cite[Proposition 5.11]{Fre}). 
\begin{proposition}\label{Frattini et Fitting}
 Let $\gk$ be a connected soluble ranked Lie ring such that $\gk'$ is nilpotent. Then $F(\gk)^{\circ}/\Phi(\gk)^{\circ}=F^{\circ}(\gk/\Phi^{\circ}(\gk))$   
\end{proposition}
\begin{proof}
Let $F$ be $\pi^{-1}(F^{\circ}(\gk/\Phi^{\circ}(\gk))$. It suffices to prove that $F$ is nilpotent. By Corollary \ref{Cartan subrings exist}, $F$ has a definable connected Cartan subring $\ck$. Since $\Phi(\gk)^{\circ}$ is nilpotent, $\pi(\ck)$ is a non-trivial Cartan subring of the nilpotent ring $F^{\circ}(\gk/\Phi^{\circ}(\gk))$, hence $F=\Phi^{\circ}(\gk)+\ck$. By Fact \ref{Frattini argument}, we have $\gk=F+N_{\gk}(\ck)=\Phi^{\circ}(\gk)+N_{\gk}(\ck)$; so $\ck$ is an ideal of $\gk$ and $F=\ck$ is nilpotent.    
\end{proof}

\section{Engel subrings}
In this section, we characterize Cartan subrings as \textit{Engel subrings}. Notice that in the classical linear finite-dimensional case, the existence of Cartan subalgebras is proved by considering minimal Engel subalgebras (see \cite{Bar}). We also take inspiration from Frécon's strategy used in the case of soluble groups of finite Morley rank (see \cite{Fre}).
\begin{definition}
Let $\gk$ be a ranked Lie ring. For all $x \in \gk$, we define the Engel subgroup associated with $x$ by \[E_{\gk}(x) = \bigcup\limits_{n\in\N^*}(\ad_x^n)^{-1}(0).\]
and for $X\subseteq \gk$, we let \[E_{\gk}(X)=\bigcap\limits_{x\in X}E_{\gk}(x).\]
\end{definition}
\begin{fait}
Let $\gk$ be a Lie ring. Then $E_{\gk}(x)$ is a subring.
\end{fait}
\begin{proof}
Let $y,z\in E_\gk (x)$. Then there exists $n,m\in \mathbb{N}$ such that $\ad^{n}_{x}(y)=\ad^{m}_x(z)=0$ and $\ad^{n-1}_{x}(y), \ad^{m-1}_x(z)\neq 0$. Set $s=n+m$. By the Leibniz product rule, we get :
\[\ad_x^s([y,z])=\sum_{0\leq k \leq s}\binom{s}{k}[\ad^{s-k}_x(y),\ad^k_x(z)]=0.\]
\end{proof}
Firs, we prove a characterization of the Fitting ideal.
\begin{lem}\label{Fitting is the set of ad-nilpotent}
Let $\gk$ be a connected soluble ranked Lie ring, such that $\operatorname{char}(\gk) >\rk(\gk)$. Then the Fitting ideal $F(\gk)$ is the set of ad-nilpotent elements of $\gk$. 
\end{lem}

\begin{proof}
It is clear that every element of $F(\gk)$ is $\ad$-nilpotent. Let $V$ be a $\gk$-minimal module on which $\gk$ acts non-trivially. Linearizing the action (Fact \ref{Lie theorem}), there is a definable field $\K$ such that $V\cong \K_+$ definably, and $\gk/C_{\gk}(V)$ acts by scalar operators. 

Let $x$ be an ad-nilpotent element of $\gk$ and $\dk(x)$ the smallest definable subring containing $x$. Then $x$ is nilpotent on $V$ and by the scalar action, $x$ centralizes $V$, and so does $\dk(x)$. Now consider a composition series 
\[0=\ak_1<\ak_2<\cdots<\ak_n=F^{\circ}(\gk)\] of $F^{\circ}(\gk)$. By the previous, we have $[\dk(x),_nF^{\circ}(\gk)]=0$. Therefore $\dk(x)+F^{\circ}(\gk)$ is nilpotent by \cite[Proposition 4.7]{Hall}, so contained in $F(\gk)$ since $\gk'$ is. We can conclude.
\end{proof}

The assumption on the characteristic at this level is important because there exist soluble Lie algebras defined over a field of positive characteristic and such that the Fitting does not contain all ad-nilpotent elements. (see \cite{Hart}, and \cite{Zam1} for an explicit example)

\begin{proposition}\label{generalized centralizers are def}
Let $\gk$ be a connected soluble ranked Lie ring, such that $\operatorname{char}(\gk)> \rk(\gk)$. Let $\hk$ be a definable connected nilpotent subring. Then $E_{\gk}(\hk)$ is definable and connected and there is $k\in\N$ such that $\ad_x^k(E_{\gk}(\hk))=0$ for all $x\in \hk$.
\end{proposition}

\begin{proof}
By induction on $\rk(\gk)$. Let $\hk$ be a definable connected nilpotent subring. We may assume that $\hk$ is not contained in $F(\gk)^{\circ}$. Suppose $\Phi^{\circ}(\gk)\neq0$, then applying induction to $\overline{\gk}=\gk/\Phi^{\circ}(\gk)$, $E_{\overline{\gk}}(\overline{\hk})$ is definable and connected and there exists $n\in \mathbb{N}$ such that $\ad^n_{\ov{x}}(E_{\overline{\gk}}(\overline{\hk}))=\ov{0}$ for all $\ov{h}\in \ov{\hk}$. If $\overline{\gk}=E_{\overline{\gk}}(\overline{\hk})$, then $\overline{\hk}\subseteq F^{\circ}(\overline{\gk})$ by Lemma  \ref{Fitting is the set of ad-nilpotent}. Since $F(\gk)^{\circ}/\Phi(\gk)^{\circ}=F^{\circ}(\gk/\Phi^{\circ}(\gk))$ (Proposition \ref{Frattini et Fitting}), $\hk\subseteq F^{\circ}(\gk)$; a contradiction. Thus $E_{\overline{\gk}}(\ov{\hk})$ is proper and there is $k\in\N$ such that 
\[\ak=\bigcap_{x\in\hk}\bigcup\limits_{n\geq0}\ad_x^{-n}(\Phi^{\circ}(\gk))=\bigcap_{x\in\hk}\ad_x^{-k}(\Phi^{\circ}(\gk))\] is a definable proper subring of $\gk$ which contains $E_{\gk}(\hk)$. By induction we can conclude that $E_{\gk}(\hk)=E_{\ak}(\hk)$ is definable and connected and there exists $k\in \mathbb{N}$ such that $\ad_h^k(\gk)=0$ for all $h\in \hk$. Now we can assume that $\Phi^{\circ}(\gk)=0$. We can also assume that $Z^{\circ}(\gk)=0$, otherwise, by applying induction to $\gk/Z^{\circ}(\gk)$, it comes that $\gk/Z^{\circ}(\gk)=E_{\gk/Z^{\circ}(\gk)}(\hk+Z(\gk)^{\circ}/Z(\gk)^{\circ})$ and there exists $n\in \mathbb{N}$ such that $\ad_h^n(\gk)\subseteq Z(\gk)^{\circ}$ for all $h\in \hk$, hence we are done.

Let $\ik$ be a $\gk$-minimal ideal contained in $\gk'$. Since $\Phi^{\circ}(\gk)=0$, we may consider $\mk$ a definable connected proper subring maximal as such which does not contain $\ik$. By Fact \ref{irreducibility of minimal modules}, $\ik$ is $\mk$-irreducible since it is not central, and we have $\gk=\ik\oplus \mk$. If $\ik< C_{\gk}^{\circ}(\ik)$, then $\jk=\mk\cap C_{\gk}^{\circ}(\ik)$ is an infinite ideal. Applying induction to $\gk/\ik$ and $\gk/\jk$, it comes that $E_{\gk/\ik}(\hk+\ik/\ik)$ and $E_{\gk/\jk}(\hk+\jk/\jk)$ are definable and connected, thus $\ak=\pi^{-1}(E_{\gk/\ik}(\hk+\ik/\ik))$ and $\bk=\pi^{-1}(E_{\gk/\jk}(\hk+\jk/\jk))$ are definable and connected and both contain $E_{\gk}(\hk)$. If one of the them is proper, induction applies and the result follows. We can suppose that $\ak=\gk=\bk$. Then, there are $k,s\in\N$ such that $\ad_h^k(\gk)\subseteq \ik$ and $\ad_h^s(\gk)\subseteq \jk$ for all $h\in \hk$. Then $\ad_h^{\max(k,s)}(\gk)\in \ik\cap \jk=\{0\}$ for all $h\in \hk$ and we are done. 

Therefore, we can assume that $\ik=C_{\gk}^{\circ}(\ik)$. As $\gk'$ is nilpotent, we have $\gk'\subseteq C_{\gk}(\ik)^{\circ}=\ik$ and $\mk$ is abelian.
By induction, we may suppose that $E_{\gk/\ik}(\hk+\ik/\ik)=\gk/\ik$. Note that we get $\gk=E_{\gk}(\hk)$ if $\hk$ centralizes $\ik$. So assume that exists $h\in\hk$ that does not centralize $\ik$. Then by Fact \ref{Lie theorem} we have $C_\ik(h)=0$. Now write $h=i+m$ with $i\in \ik$ and $m\in\mk$. We have $C_\ik(m)=0$ and the map $f(z)=[z,m]$ is injective on $\ik$ thus surjective and we can write $i=[z,m]$. This implies $h\in \mk_1$ where $\mk_1=\exp(\ad_z)(\mk)$. Then $\mk_1\sqsubseteq E_{\gk}(h)<\gk$ since $\ik\cap \mk_1=0$ and $E_{\gk/\ik}(\hk+\ik/\ik)=\gk/\ik$. In fact, we get $\mk_1=E_{\gk}(h)$ : let $x\in E_{\gk}(h)$; we have $x=i+m$ for some $i\in \ik$ and $m\in\mk_1$, and so there exists $k\in \mathbb{N}$ such that $\ad^k_h(i+m)=\ad_h^k(i)=0$ and $i=0$ since $C_{\ik}(h)=0$. In particular $\hk \sqsubseteq E_{\gk}(\hk)\sqsubseteq E_{\gk}(h)=\mk_1$. But for all $h\in\hk$ and $m\in\mk_1$, there is $n\in\N$ such that $\ad_h^n(m)\in \ik\cap \mk_1=0$ since $E_{\gk/\ik}(\hk+\ik/\ik)=\gk/\ik$; so $\mk_1=E_{\gk}(\hk)$. Finally, note that $C_{\gk}(\hk)=\mk=E_{\gk}(\hk)$ since $\mk$ is abelian and this ends the proof.
\end{proof}

\begin{definition}
Let $V$ be a $\gk$-module. We say that $x\in \gk$ is \emph{$V$-nilpotent} if there is $n\in\N^*$ such that $x^n\cdot V=0$. We say that $\gk$ is $V$-nilpotent if every $x\in \gk$ is $V$-nilpotent.
\end{definition}

In the literature, it is also said that $\gk$ is locally nilpotent on $V$. In the case of the derivation action on $\gk$, $\gk$-nilpotency means $\ad$-nilpotency.

\begin{lem}\label{nilpotent trivial action on irreducuble}
Let $\gk$ be a connected nilpotent ranked Lie ring and let $V$ be a connected definable $\gk$-module of characteristic $p>0$. Suppose that $p > \rk(V)$ and that $\gk$ is $V$-nilpotent. Then $V$ contains a non-trivial $\gk$-submodule on which $\gk$ acts trivially.
\end{lem}

\begin{proof}
We can assume that $V$ is $\gk$-minimal. By Fact \ref{Lie theorem}, there is a definable field $\K$ such that $V\cong\K_+$ definably and $\gk/C_{\gk}(V)$ acts by scalar operators. Since $\gk$ is $V$-nilpotent, it follows immediately that $\gk$ acts trivially on $V$.
\end{proof}

We can now prove, in the following, that Cartan subrings in our context are Engel subrings.

\begin{proposition}\label{Cartan subring is Engel}
Let $\gk$ be a connected soluble ranked Lie ring such that $\operatorname{char}(\gk) > \rk(\gk)$. Let $\hk$ be a Cartan subring. Then $E_{\gk}(\hk) = \hk$.
\end{proposition}

\begin{proof}
Suppose that the $\hk$-module $V=E_{\gk}(\hk)/\hk$ is not trivial. By Proposition \ref{generalized centralizers are def}, $\hk$ is $V$-nilpotent, then by Lemma \ref{nilpotent trivial action on irreducuble}, there is a non-trivial $\hk$-submodule $W=\ak/\hk$ on which $\hk$ acts trivially. This implies that $[\ak,\hk]\subseteq\hk$, so $\ak\sqsubseteq \hk=N_{\gk}(\hk)$; this is absurd since $W$ is non-trivial. We can conclude.
\end{proof}

\begin{proposition}\label{Engel's set lifts}
Let $\gk$ be a connected soluble ranked Lie ring such that $\operatorname{char}(\gk) > \rk(\gk)$. Let $\ik$ be a definable ideal and $\hk$ be a definable connected nilpotent subring. Then $E_{\gk/\ik}(\hk + \ik/\ik) = E_{\gk}(\hk) + \ik/\ik$.
\end{proposition}

\begin{proof}
We proceed by induction on $\rk(\gk)$. We can assume that $\gk$ is not nilpotent. We can also assume that $\ik$ is $\gk$-minimal (hence connected) and not central. If $E_{\gk/\ik}(\hk+\ik/\ik)$ is proper, then applying induction to its preimage we are done. Therefore, we can suppose that $E_{\gk/\ik}(\hk+\ik/\ik)=\gk/\ik$. Finally we can assume that $\ik$ is not contained in $\Phi(\gk)^{\circ}$. Indeed, if $E_{\gk/\Phi(\gk)^{\circ}}(\hk+\Phi(\gk)^{\circ}/\Phi(\gk)^{\circ})=\gk/\Phi(\gk)^{\circ}$, by Proposition \ref{Frattini et Fitting}, $\hk\subseteq F(\gk)^{\circ}$; hence we are done by applying induction to $E_{\gk/\Phi(\gk)^{\circ}}(\hk+\Phi(\gk)^{\circ}/\Phi(\gk)^{\circ})$.

Let $\mk$ be a connected maximal subring not containing $\ik$. Then $\gk=\ik\oplus \mk$ and $\ik$ is $\mk$-irreducible by Lemma \ref{irreducibility of minimal modules} since we may suppose that $\ik$ is not central. Note that we get $\gk=E_{\gk}(\hk)$ if $\hk$ centralizes $\ik$. So assume that exists $h\in\hk$ that does not centralize $\ik$. Then by Fact \ref{Lie theorem} we have $C_\ik(h)=0$. Now write $h=i+m$ with $i\in \ik$ and $m\in\mk$. We have $C_\ik(m)=0$ and the map $f(z)=[z,m]$ is injective on $\ik$ thus surjective and we can write $i=[z,m]$. This implies $h\in \mk_1$ where $\mk_1=\exp(\ad_z)(\mk)$. Then $\mk_1\sqsubseteq E_{\gk}(h)<\gk$ since $\ik\cap \mk_1=0$ and $E_{\gk/\ik}(\hk+\ik/\ik)=\gk/\ik$. In fact, we get $\mk_1=E_{\gk}(h)$ : let $x\in E_{\gk}(h)$; we have $x=i+m$ for some $i\in \ik$ and $m\in\mk_1$, and so there exists $k\in \mathbb{N}$ such that $\ad^k_h(i+m)=\ad_h^k(i)=0$ and $i=0$ since $C_{\ik}(h)=0$. In particular $\hk\sqsubseteq E_{\gk}(\hk)\sqsubseteq E_{\gk}(h)\sqsubseteq\mk_1$. But for all $h\in\hk$ and $m\in\mk_1$, there is $n\in\N$ such that $\ad_h^n(m)\in \ik\cap \mk_1=0$ since $E_{\gk/\ik}(\hk+\ik/\ik)=\gk/\ik$; so $\mk_1=E_{\gk}(\hk)$. Therefore we have $\gk=\ik\oplus E_{\gk}(\hk)$ and this ends the proof.
\end{proof}

So far, we know that a minimal def-abnormal subring is Cartan (Proposition \ref{minimal abnorma is Cartan}) and Cartan subrings are Engel subrings. We finally prove that Engel subrings are def-abnormal subrings. 

\begin{proposition}\label{Engel subring is abnormal}
Let $\gk$ be a connected soluble ranked Lie ring such that $\operatorname{char}(\gk) > \rk(\gk)$. Let $\hk$ be a definable connected nilpotent subring. Then $E_{\gk}(\hk)$ is def-abnormal in $\gk$.
\end{proposition}

\begin{proof}
 We proceed by induction on $\rk(\gk)$. Let $\ik$ be a $\gk$-minimal ideal. By induction and Proposition \ref{Engel's set lifts} $E_{\gk/\ik}(\hk+\ik/\ik)=E_{\gk}(\hk)+\ik/\ik$ is def-abnormal in $\gk/\ik$, so $E_{\gk}(\hk)+\ik$ is def-abnormal in $\gk$. By Proposition \ref{abnormality is transitive} it suffices to show that $E_{\gk}(\hk)$ is def-abnormal in $E_{\gk}(\hk)+\ik$. But $\gk/C_{\gk}(\ik)$ is abelian by Fact \ref{Lie theorem}, so 
\[\gk/C_{\gk}(\ik)=E_{\gk/C_{\gk}(\ik)}(\hk+C_{\gk}(\ik)/C_{\gk}(\ik))=E_{\gk}(\hk)+C_{\gk}(\ik)/C_{\gk}(\ik).\] Thus $\gk=E_{\gk}(\hk)+C_{\gk}(\ik)$ and we can conclude by Proposition \ref{abnormality criterion}.
\end{proof}

\begin{definition}
A subring $\ak$ is said to be \emph{Engel-minimal} if there exists a definable connected nilpotent subring $\hk$ such that $\ak = E_{\gk}(\hk)$ and if $\ak$ is minimal among subrings of the form $E_{\gk}(\mathfrak{n})$ where $\mathfrak{n}$ is a definable connected nilpotent subring.
\end{definition}

The following ends up establishing the equivalence between Cartan subrings, minimal def-abnormal subrings and Engel-minimal subrings.

\begin{proposition}\label{Engel-minimal is Cartan}
Let $\gk$ be a connected soluble ranked Lie ring such that $\operatorname{char}(\gk) > \rk(\gk)$. Let $\ak$ be an Engel-minimal subring. Then $\ak$ is a Cartan subring.
\end{proposition}

\begin{proof}
Note that $\ak$ is def-abnormal by Proposition \ref{Engel subring is abnormal}. By Proposition \ref{abnormal exits}, $\ak$ has a minimal def-abnormal subring $\ck$ which is a Cartan subring (Corollary \ref{minimal abnorma is Cartan}). In particular $\ck$ is def-abnormal in $\gk$ by transitivity (Proposition \ref{abnormality is transitive}). But $\ck=E_{\gk}(\ck)$ (Proposition \ref{Cartan subring is Engel}), therefore $\ck=\ak$ by minimality of $\ak$. 
\end{proof}

We finally summarize all this in the following theorem.

\begin{theo}\label{main theorem}
Let $\gk$ be a connected soluble ranked Lie ring such that $\operatorname{char}(\gk)>\rk(\gk)$ and let $\ck\sqsubseteq\gk$ be a subring. The following items are equivalent:
\begin{itemize}
    \item [(i)] $\ck$ is a minimal def-abnormal subring;
    \item [(ii)] $\ck$ is a Cartan subring;
    \item [(iii)] $\ck$ is an Engel-minimal subring.
\end{itemize}
\end{theo}

\begin{proof}
    First note that $(i)\Rightarrow (ii)$ is Proposition \ref{minimal abnorma is Cartan}; $(iii)\Rightarrow (ii)$ is Proposition \ref{Engel-minimal is Cartan}; and $(ii)\Rightarrow(iii)$ is Proposition \ref{Cartan subring is Engel} together with Proposition \ref{Engel subring is abnormal}. 
   
    Let $\ak$ be an Engel-minimal subring. Then $\ak$ is def-abnormal  (Proposition \ref{Engel subring is abnormal}) and if $\ck$ is a minimal def-anormal subring contained in $\ak$, then $\ck$ is Cartan, hence $\ck=E_{\gk}(\ck)$ (Proposition \ref{Cartan subring is Engel}). Therefore $\ak=\ck$ by minimality, so $\ak$ is a minimal def-abnormal subring and $(iii)\Rightarrow (i)$ is done.
\end{proof}
\bibliography{bib}
\bibliographystyle{plain}
\end{document}